\documentclass{amsart}

\newtheorem{theorem}{Theorem}
\newcommand{\bt}{\begin{theorem}}
\newcommand{\et}{\end{theorem}}

\newcommand{\beq}{\begin{equation}}
\newcommand{\eeq}{\end{equation}}
\newcommand{\benum}{\begin{enumerate}}
\newcommand{\eenum}{\end{enumerate}}

\newcommand{\Z}{\ensuremath{\mathbf Z}}

\newcommand{\N}{\ensuremath{ \mathbf N }}
\newcommand{\mci}{\ensuremath{\mathcal I}}

\newcommand{\mcz}{\ensuremath{\mathcal Z}}

\begin{document}

\title{Adjoining identities and zeros to semigroups}

\author{Melvyn B. Nathanson}
\address{Lehman College (CUNY),Bronx, New York 10468}
\email{melvyn.nathanson@lehman.cuny.edu}

\thanks{This work was supported in part by grants from the NSA Mathematical Sciences Program and the PSC-CUNY Research Award Program.}

\date{\today}

\maketitle

\begin{abstract}
This note shows how iteration of the standard process of adjoining identities and zeros to semigroups gives rise naturally to the lexicographical ordering on the additive semigroups of $n$-tuples of nonnegative integers and $n$-tuples of integers.
\end{abstract}

\section{Semigroups with identities and zeros}
A binary operation $\ast$ on a set $S$ is \emph{associative} if $(a\ast b)\ast c = a \ast (b\ast c)$ for all $a,b,c \in S$.  A \emph{semigroup} is a nonempty set  with an associative binary operation $\ast$.  The semigroup is \emph{abelian} if $a\ast b = b\ast a$ for all $a,b \in S.$  The trivial semigroup $S_0$ consists of a single element $s_0$ such that $s_0 \ast s_0 = s_0$.  Theorems about abstract semigroups are, in a sense, theorems about the pure process of multiplication.

An element $u$ in a semigroup $S$ is an \emph{identity} if $u\ast a = a\ast u = a$ for all $a \in S.$  If $u$ and $u'$ are identities in a semigroup, then $u = u\ast u' = u'$ and so a semigroup contains at most one identity.  A semigroup with an identity is called a \emph{monoid}.
If $S$ is a semigroup that is not a monoid, that is, if $S$ does not contain an identity element, there is a simple process to adjoin an identity to $S$.  Let $u$ be an element not in $S$ and let 
\[
\mci(S,u) = S \cup \{u\}.
\]
We extend the binary operation $\ast$ from $S$ to $\mci(S,u)$ by defining $u\ast a = a \ast u = a$ for all $a \in S,$ and $u\ast u = u.$  Then $\mci(S,u)$ is a monoid with identity $u$.

An element $v$ in a semigroup $S$ is a \emph{zero} if $v\ast a = a\ast v = v$ for all $a \in S.$  If $v$ and $v'$ are zeros in a semigroup, then $v=v\ast v'=v'$ and so a semigroup contains at most one zero.  If $S$ is a semigroup that does not contain a zero element, there is also a simple process to adjoin a zero to $S$.  Let $v$ be an element not in $S$ and let 
\[
\mcz(S,v) = S \cup \{v\}.
\]
We extend the binary operation $\ast$ from $S$ to $\mcz(S,v)$ by defining $v\ast a = a \ast v = v$ for all $a \in S$, and $v\ast v = v.$  Then $\mcz(S,v)$ is a semigroup with zero $v$.

It is important to note that the process of adjoining an identity to a semigroup $S$ is well-defined even if $S$ contains an identity.  Similarly, the process of adjoining a zero to a semigroup $S$ is well-defined even if $S$ contains a zero.  

The element $s_0$ in the trivial semigroup $S_0$  is both an identity and a zero.  In this note we investigate what happens when we start with the trivial semigroup and add new identities and new zeros.

\section{Lexicographical order}
A relation $\leq$ on a set $S$ is a \emph{partial order} if 
\benum
\item[(A1)]
$a \leq a$ for all $a\in S$
\item[(A2)]
If $a,b,c \in S$ and $a\leq b$ and $b\leq c$, then $a\leq c$
\item[(A3)]
If $a,b \in S$ and $a\leq b$ and $b\leq a$, then $a=b$.
\eenum
The relation is called a \emph{total order} if it satisfies the additional  condition:
\benum
\item[(A4)]
If $a,b \in S$, then $a\leq b$ or $b\leq a$.
\eenum
We write $a<b$ if $a \leq b$ and $a\neq b$.

Let $S$ be a partially ordered set and let $\{X_s\}_{s\in S}$ be a family of pairwise disjoint nonempty partially ordered sets.  We define a partial order on the set $X = \cup_{s\in S} X_s$ as follows:  Let $x,y \in X$.  We write $x \leq y$ if  there exists $s \in S$ such that $a,b \in X_s$ and $a \leq b$ in $X_s$, or if there exist $s_1, s_2 \in S$ with $s_1 < s_2$ such that  $a\in X_{s_1}$ and $b \in X_{s_2}$.  If $S$ is a totally ordered set and if $X_s$ is totally ordered for each $s\in S$, then the set $X$ is also totally ordered.
We have constructed $X$ by replacing each element of $S$ with the set $X_s.$

Let $\leq$ be a total order on $S$.  If $a,b \in S$ and $a \leq b$,  we define $\min(a,b) = a$.
We define a binary operation $a\ast b$ on $S$ by $a\ast b = \min(a,b)$.  This operation has the following properties:
\benum
\item[(B1)]  $(a \ast b) \ast c = a \ast (b \ast c)$ and $a \ast b = b \ast a$ for all $a,b,c \in S$, and so $S$ is an abelian semigroup with respect to $\ast$.
\item[(B2)]
$a\ast b \in \{a,b\}$ for all $a,b\in S$.
\item[(B3)]
If $a\ast b = a$ and $b\ast c = b$, then $a \ast c = a$.
\eenum
The binary operation on $S$ defined by $a\ast b = \max(a,b)$ also satisfies~(B1),~(B2), and~(B3).  
Conversely, if $S$ is a semigroup whose binary operation $\ast$ satisfies properties~(B1),~(B2), and~(B3), and if we define a relation $\leq$ on $S$ by
$a \leq b$ if $a\ast b = a$, then $\leq$ is a total order on $S$.

Let  $\N = \{1,2,3\ldots\}$ denote the set of positive integers, $\N_0 = \N \cup \{0\} $  the set of nonnegative integers, and \Z\ the set of all integers.
For every positive integer $n$, we let $\N_0^n$ and $\Z^n$ denote the sets of $n$-tuples of nonnegative integers and integers, respectively.  These are abelian semigroups with respect to the operation $+$ of ordinary vector addition.  

We define the \emph{lexicographic order} on $\Z^n$ as follows.   Let $(i_1,\ldots, i_n)\in \Z_0^n$ and $(j_1,\ldots, j_n)\in \Z_0^n$.  We define $(i_1,\ldots, i_n) \leq (j_1,\ldots, j_n)$ if either  $(i_1,\ldots, i_n) = (j_1,\ldots, j_n)$ or if  $(i_1,\ldots, i_n) \neq (j_1,\ldots, j_n)$ and $i_k < j_k$, where $k$ is the smallest positive integer such that $i_k \neq j_k.$
Lexicographic order is a total order on $\Z^n$ and, by restriction, on $\N_0^n$.

A total order $\leq$ on $\N_0^n$ is called a \emph{monomial order} if it is a well-order, that is, every nonempty subset of $\N_0^n$ contains a smallest element, and if, for all $a,b,c\in \N_0^n$ with $a \leq b$, we have $a+c \leq b+c$.  The lexicographical order is a monomial order on $\N_0^n$ (Cox, Little, O'Shea~\cite[Section 2.2]{cox-litt-oshe97}).

The object of this note is to show that, starting with the trivial semigroup and iterating the process of adjoining an identity, we obtain the set $\N_0^n$ with the lexicographical ordering.  If we adjoin infinitely many identities and infinitely many zeros to the trivial semigroup, we obtain the set $\Z$ with the usual ordering.  Replacing each element of \Z\ with a copy of \Z\ generates $\Z^2$ with the lexicographical ordering.  Iterating this process finitely many times produces $\Z^n$  with the lexicographical ordering.

\section{Iteration of the process of adjoining an identity}

Let $S_0 = \{s_0\}$ be the trivial semigroup.  For $s_1 \neq s_0$, let $S_1 = \mci(S_0,s_1) = \{s_0,s_1\}$ be the semigroup obtained by adjoining the identity $s_1$ to $S_0$.  Then $s_0\ast s_1 = s_1 \ast s_0 = s_0$ and $s_1 \ast s_1 = s_1.$   For $s_2 \neq s_0, s_1$, let $S_2 = \mci(S_1,s_2) = \{s_0,s_1,s_2\}$ be the semigroup obtained by adjoining the identity $s_2$ to $S_1$.  
Then $s_i\ast s_2 = s_2 \ast s_i = s_i$ for $i=0,1,2.$ 
Continuing inductively, we obtain an increasing sequence of abelian semigroups
\[
S_0 \subseteq S_1 \subseteq S_2 \subseteq \cdots \subseteq S_k \subseteq \cdots
\]
such that 
\[
S_k = \{s_0, s_1, s_2,\ldots, s_k\}
\]
for all $k\in \N_0$.  Then 
\[
T^{(1)} = \bigcup_{k=0}^{\infty}S_k = \{s_0,s_1,s_2,\ldots,s_i,\ldots\}
\]
is a abelian semigroup with the binary operation 
\[
s_i\ast s_{j} = s_{\min(i,j)}
\]
for all $i,j \in \N_0.$ 
Note that $s_0$ is a zero in the semigroup $T^{(1)}$, but that this semigroup does not contain an identity.

Define $s_{0,i} = s_i$ for all $i \in \N_0$ and let 
\[
U^{(0)} = \{s_{0,0}, s_{0,1}, s_{0,2}, \ldots, s_{0,i},\ldots\}.
\]
Choose $s_{1,0} \notin U^{(0)}$ and consider the semigroup
\[
U^{(0)}_0 = \mci(U^{(0)},s_{1,0}) = \{s_{0,0}, s_{0,1}, s_{0,2}, \ldots, s_{0,i},\ldots,s_{1,0}\}.
\]
For all $i,j \in \N_0$ we have
\[
s_{0,i} \ast s_{0,j} = s_{0,\min(i,j)}
\]
and
\[
s_{0,i}\ast s_{1,0} = s_{0,i}.
\]
Choose $s_{1,1} \notin U^{(0)}_0$ and let 
\[
U^{(0)}_1 = \mci(U^{(0)}_0,s_{1,1}) = \{s_{0,0}, s_{0,1}, s_{0,2}, \ldots, s_{0,i},\ldots,s_{1,0},s_{1,1}\}.
\]
Iterating this process, we obtain an increasing sequence of abelian semigroups
\[
U^{(0)}_0 \subseteq U^{(0)}_1 \subseteq U^{(0)}_2 \subseteq \cdots \subseteq U^{(0)}_k \subseteq \cdots
\]
with 
\[
U^{(0)}_k  = \{s_{0,0}, s_{0,1}, s_{0,2}, \ldots, s_{0,i}, \ldots,s_{1,0},s_{1,1},s_{1,2},\ldots,s_{1,k}\}.
\]
Then
\[
U^{(1)} = \bigcup_{k=0}^{\infty} U^{(0)}_k = \{s_{0,0}, s_{0,1}, s_{0,2}, \ldots, s_{0,i},\ldots,s_{1,0},s_{1,1},s_{1,2},\ldots,s_{1,i},\ldots\}
\]
is a abelian semigroup with the multiplication 
\[
s_{0,i} \ast s_{0,j} = s_{0,\min(i,j)}
\]
\[
s_{1,i} \ast s_{1,j} = s_{1,\min(i,j)}
\]
and
\[
s_{0,i} \ast s_{1,j} = s_{0,i}
\]
for all $i,j \in \N_0.$
Iterating this process, we obtain an increasing sequence of abelian semigroups
\[
U^{(0)} \subseteq U^{(1)} \subseteq U^{(2)} \subseteq \cdots \subseteq U^{(\ell)} \subseteq \cdots
\]
with 
\[
U^{(\ell)}  = \bigcup_{i_1=0}^{\ell}  \{s_{i_1,i_2}\}_{i_2=0}^{\infty}
\]
such that 
\[
T^{(2)} = \bigcup_{\ell =0}^{\infty} U^{(\ell)} = \{ s_{i_1,i_2} \}_{i_1,i_2=0}^{\infty}
\]
is an abelian semigroup whose multiplication satisfies
\[
s_{i_1,i_2} \ast s_{j_1,j_2} = 
\begin{cases}
s_{i_1,\min(i_2,j_2)} & \text{if $i_1=j_1$}\\
s_{i_1,i_2} & \text{if $i_1 < j_1$.}
\end{cases}
\]

Again, iterating the process of adjoining identities to semigroups, we obtain, for every positive integer $n$, the abelian semigroup 
\[
T^{(n)} = \{s_{i_1,i_2,\ldots,i_n} \}_{i_1,\ldots,i_n=0}^{\infty}
\]
with multiplication $\ast$ defined by
\[
s_{i_1,i_2,\ldots,i_n} \ast s_{j_1,j_2,\ldots,j_n}
=  
s_{i_1,i_2,\ldots,i_n} 
\]
if either  $(i_1,\ldots, i_n) = (j_1,\ldots, j_n)$ or if  $(i_1,\ldots, i_n) \neq (j_1,\ldots, j_n)$ and $i_k < j_k$, where $k$ is the smallest positive integer such that $i_k \neq j_k.$  The binary operation on $T^{(n)}$ satisfies properties~(B1),~(B2), and~(B3), and so induces a total order on this semigroup.  Since 
\[
s_{i_1,i_2,\ldots,i_n} \ast s_{j_1,j_2,\ldots,j_n}
=  s_{i_1,i_2,\ldots,i_n}
\]
if and only if 
\[
(i_1,i_2,\ldots,i_n) \leq (j_1,j_2,\ldots,j_n)
\]
with respect to the lexicographical order, it follows that the process of iterated adjunction of an identity to the trivial semigroup has recreated the semigroup $\N_0^n$ with the lexicographical order.

\section{Iteration of the process of adjoining a zero}

We return to the semigroup $T^{(1)} = \{s_0,s_1,s_2,\ldots,\}$ with the binary operation $s_i\ast s_j = s_{\min(i,j)}.$  Choose an element $s_{-1}$ such that $s_{-1} \neq s_i$ for all $i\in \N_0$ and let $T^{(1)}_1 = \mcz(T^{(1)},s_{-1})$ be the semigroup obtained by adjoining the zero $s_{-1}$ to $T^{(1)}.$  Then $s_i\ast s_{-1} = s_{-1} \ast s_i = s_{-1}$ for all integers $i\geq -1.$  Choose an element $s_{-2}$ such that $s_{-2} \neq s_i$ for all $i\in \N_0\cup \{-1\}$ and let $T^{(1)}_2 = \mcz(T^{(1)}_1,s_{-2})$ be the semigroup obtained by adjoining the zero $s_{-2}$ to $T^{(1)}_1.$  Then $s_i\ast s_{-2} = s_{-2} \ast s_i = s_{-2}$ for all integers $i\geq -2.$
Continuing inductively, we obtain an increasing sequence of abelian semigroups
\[
T^{(1)} \subseteq T^{(1)}_1 \subseteq T^{(1)}_2 \subseteq \cdots \subseteq T^{(1)}_k \subseteq \cdots
\]
such that 
\[
T^{(1)}_k = \{s_{-k},s_{-k+1},\ldots,s_{-1},s_0, s_1, s_2,\ldots, s_i,\ldots\}
\]
for all $k\in \N_0$.  Then
\[
V^{(1)} = \bigcup_{k=0}^{\infty}T^{(1)}_k = \{s_i\}_{i=-\infty}^{\infty}
\]
is a abelian semigroup with the binary operation 
\[
s_i\ast s_{j} = s_{\min(i,j)}
\]
for all $i,j \in \Z.$  
This operation satisfies properties~(B1)-(B3) and so induces a total order on $V^{(1)}.$
The semigroup $V^{(1)}$ contains neither an identity nor a zero.

Let $s_{i_1,i_2} = s_{i_2}$ for all $i_1,i_2 \in \Z$ and consider the set
\[
X_{s_{i_1}} = \{s_{i_1,i_2}\}_{i_2=-\infty}^{\infty}
\]
with the total order defined by $s_{i_1,i_2} \leq s_{i_1,j_2}$ if $i_2 \leq j_2.$
Then $\{X_{s_{i_1}}\}_{s_{i_1} \in V^{(1)}}$ is a family of pairwise disjoint nonempty totally ordered sets.   Replacing each element $s_{i_1}\in V^{(1)}$  with the  set $X_{s_{i_1}}$, we obtain the totally ordered set
\[
V^{(2)} = \bigcup_{s_{i_1} \in V^{(1)}}  X_{s_{i_1}}   = \{s_{i_1,i_2}\}_{i_1,i_2=-\infty}^{\infty}
\]
where  $s_{i_1,i_2} \leq s_{j_1,j_2}$ in $V^{(2)}$ 
if either $i_1 < j_1$ or $i_1 = j_1$ and $i_2 \leq j_2$.

Defining $s_{i_1,i_2,i_3} = s_{i_3}$ for all $i_1,i_2,i_3 \in \Z$ and replacing each element 
$s_{i_1,i_2}$ of $V^{(2)}$ with the set $X_{s_{i_1,i_2}} = \{s_{i_1,i_2,i_3}\}_{i_3=-\infty}^{\infty}$,
we obtain the totally ordered set
\[
V^{(3)} = \bigcup_{s_{i_1,i_2} \in V^{(2)} } X_{s_{i_1,i_2}} 
= \{s_{i_1,i_2,i_3}\}_{i_1,i_2,i_3=-\infty}^{\infty}
\]
with the lexicograhical ordering.
Iterating this process, we obtain, for every positive integer $n$, the abelian semigroup 
\[
V^{(n)} = \{s_{i_1,i_2,\ldots,i_n} \}_{i_1,\ldots,i_n=-\infty}^{\infty}
\]
with multiplication $\ast$ defined by
\[
s_{i_1,i_2,\ldots,i_n} \ast s_{j_1,j_2,\ldots,j_n}
=  s_{i_1,i_2,\ldots,i_n} 
\]
if either  $(i_1,\ldots, i_n) = (j_1,\ldots, j_n)$ or if  $(i_1,\ldots, i_n) \neq (j_1,\ldots, j_n)$ and $i_k < j_k$, where $k$ is the smallest positive integer such that $i_k \neq j_k.$  The binary operation on $V^{(n)}$ satisfies properties~(B1),~(B2), and~(B3), and so induces a total order on this semigroup.  Since 
\[
s_{i_1,i_2,\ldots,i_n} \ast s_{j_1,j_2,\ldots,j_n}
=  s_{i_1,i_2,\ldots,i_n}
\]
if and only if 
\[
(i_1,i_2,\ldots,i_n) \leq (j_1,j_2,\ldots,j_n)
\]
with respect to the lexicographical order, it follows that the processes of iterated adjunction of  identities and zeros to the trivial semigroup and iterated replacement of partially ordered sets produces the semigroup $\Z^n$ with the lexicographical order.
In this way we have created something from nothing.

\emph{Remarks.}
For more information about the general theory of semigroups, there are excellent books by  Clifford and Preston~\cite{clif-pres61,clif-pres67}, Howie~\cite{howi95}, and Ljapin~\cite{ljap74}.

I wish to thank Richard Bumby for helpful discussions on this topic.

\def\cprime{$'$}
\providecommand{\bysame}{\leavevmode\hbox to3em{\hrulefill}\thinspace}
\providecommand{\MR}{\relax\ifhmode\unskip\space\fi MR }
\providecommand{\MRhref}[2]{%
  \href{http://www.ams.org/mathscinet-getitem?mr=#1}{#2}
}
\providecommand{\href}[2]{#2}


\begin{thebibliography}{1}

\bibitem{clif-pres61}
A.~H. Clifford and G.~B. Preston, \emph{The algebraic theory of semigroups.
  {V}ol. {I}}, Mathematical Surveys, No. 7, American Mathematical Society,
  Providence, R.I., 1961. \MR{MR0132791 (24 \#A2627)}

\bibitem{clif-pres67}
\bysame, \emph{The algebraic theory of semigroups. {V}ol. {II}}, Mathematical
  Surveys, No. 7, American Mathematical Society, Providence, R.I., 1967.
  \MR{MR0218472 (36 \#1558)}

\bibitem{cox-litt-oshe97}
David Cox, John Little, and Donal O'Shea, \emph{Ideals, varieties, and
  algorithms}, second ed., Undergraduate Texts in Mathematics, Springer-Verlag,
  New York, 1997, An introduction to computational algebraic geometry and
  commutative algebra. \MR{MR1417938 (97h:13024)}

\bibitem{howi95}
John~M. Howie, \emph{Fundamentals of semigroup theory}, London Mathematical
  Society Monographs. New Series, vol.~12, The Clarendon Press Oxford
  University Press, New York, 1995, , Oxford Science Publications.
  \MR{MR1455373 (98e:20059)}

\bibitem{ljap74}
E.~S. Ljapin, \emph{Semigroups}, third ed., American Mathematical Society,
  Providence, R.I., 1974, Translated from the 1960 Russian original by A. A.
  Brown, J. M. Danskin, D. Foley, S. H. Gould, E. Hewitt, S. A. Walker and J.
  A. Zilber, Translations of Mathematical Monographs, Vol. 3. \MR{MR0352302 (50
  \#4789)}

\end{thebibliography}
\end{document}